\documentclass[11pt,a4paper,reqno]{amsart}
\newif\ifchfive \chfivetrue

\newtheorem{proposition}{Proposition}

\newtheorem{theorem}{Theorem}

\newtheorem{lemma}{Lemma}

\newtheorem{conjecture}{Conjecture}

\def\eps{\varepsilon}
\def\R{{\mathbb R}}

\def\N{{\mathbb N}}

\def\O{{\mathcal O}}

\def\P{{\mathcal P}}
\def\1{ {\hbox{{\it 1}} \!\! I} }
\def\1{{{\mathit 1} \!\!\>\!\! I} }
\renewcommand{\phi}{\varphi}

\def\orb{\hbox{{\rm orb}}}

\def\ph{\varphi}

\begin{document}
\author{Henk Bruin}
\address{Department of Mathematics, 
University of Surrey, UK}
\email{{h.bruin@eim.surrey.ac.uk}}
\urladdr{{http://www.maths.surrey.ac.uk/showstaff?H.Bruin}}
 
\author{Stefano Luzzatto}
\address{Department of Mathematics,
Imperial College, London}
\email{{stefano.luzzatto@imperial.ac.uk}}
\urladdr{{http://www.ma.ic.ac.uk/\textasciitilde luzzatto}}

\thanks{We would like to thank Juan Rivera-Letelier for drawing our attention
to previous results relating to  
Lemma~\ref{equidistribution} and 
Proposition~\ref{lowlyap}. We also thank Feliks Przytycki for pointing
out an error in an earlier version of this paper.}
\date{\today}

\title[Topological invariance of the sign of Lyapunov exponents]{
Topological invariance of the sign of the Lyapunov exponents
in one-dimensional maps}

\begin{abstract}
We explore some properties of Lyapunov exponents of measures preserved by 
smooth maps of the interval, and study the behaviour 
of the Lyapunov exponents under topological conjugacy. 
\end{abstract}

\maketitle

\section{Statement of results}\label{sec_intro}

In this paper we consider \( C^3 \) interval maps \( f: I \to I \) 
where \( I \) is a compact interval. We let \( \mathcal C \) denote 
the set of \emph{critical points} of \( f \): \( c\in\mathcal C 
\Leftrightarrow Df(c) = 0 \). We shall always suppose that  
\( \mathcal C \) is finite and that each critical point
is non-flat: for each $c \in \mathcal C$, there exist 
$\ell = \ell(c) \in [2,\infty)$
and $K$ such that $\frac1K \leq \frac{|f(x)-f(c)|}{|x-c|^\ell} \leq K$
for all $x \neq c$.
Let ${\mathcal M}$ be the set of 
ergodic Borel $f$-invariant probability measures. 
For every $\mu \in \mathcal M$, we define the 
\emph{Lyapunov exponent} \( \lambda(\mu) \)  by 
 \[
\lambda(\mu) = \int \log |Df| d\mu.
\]
Notice that \( \int \log |Df| d\mu < + \infty \) is automatic 
since \( Df \) is bounded. 
 However we can have \( \int \log |Df| d\mu = - \infty \) if
\( c\in  \mathcal C \) is a fixed point and \( \mu \) is the Dirac-$\delta$ 
 measure on \( c \). 
It follows from \cite{Prz93, Blokh83} that this
is essentially the only way in which $\log |Df|$ can be non-integrable:
if $\mu(\mathcal C) = 0$, then \( \int \log |Df| d\mu > - \infty \).

 The sign, more than the actual value, of the Lyapunov exponent can 
 have significant implications for the dynamics. A positive Lyapunov 
 exponent, for example, indicates sensitivity to initial conditions 
 and thus ``chaotic'' dynamics of some kind. Our main result concerns 
 the extent to which the sign of the Lyapunov exponent, which is 
 a priori a purely \emph{metric} condition, is in fact 
intrinsically constrained by the 
\emph{topological} structure of the dynamics. 

\begin{theorem}\label{sign_lya_meas}
If $f$ is $C^3$ with finitely many non-flat critical points, and
if \( \mu \) is non-atomic then the sign of  
\( \lambda(\mu) \) is a topological invariant.
\end{theorem}

We recall that \( \mu \) is \emph{non-atomic} if every point 
has zero measure. By the statement that the sign of \( \lambda(\mu) \) 
is a topological invariant we mean the following. Two maps \( f: I \to I \) and  
\( g: J \to J \) are \emph{topologically conjugate} if there exists a 
homeomorphism \( h: I \to J \) such that \( h\circ f = g \circ h \). 
The conjugacy \( h \) induces a bijection between the space of ergodic 
invariant probability measures of \( f \) and of \( g \): if \( 
\mu_{f} \) is an ergodic invariant probability measure for \( f \), 
then the corresponding measure \( \mu_{g} \), defined 
by  \( \mu_{g}(A) = \mu_{f}(h^{-1}(A)) \) for all 
measurable  sets \( A \), is an ergodic invariant probability measure 
for \( g \). Theorem \ref{sign_lya_meas} says that as long as 
both $f$ and $g$ are $C^3$ with finitely many non-flat critical points
and \( \mu_{f} \) is non-atomic,  
then the Lyapunov exponents \( \lambda(\mu_{f}) \) and \( 
\lambda(\mu_{g}) \) have the same sign. 
Clearly the actual values can vary.

The non-atomic condition is necessary in general as a topological 
conjugacy can easily map a hyperbolic attracting/repelling periodic point to a 
topologically attracting/repelling\footnote{If $f$ has
negative Schwarzian derivative, then a neutral periodic point cannot be 
two-sided repelling} neutral periodic point. The corresponding 
Lyapunov exponents of the corresponding Dirac-$\delta$ measures 
would then be positive and zero respectively. The result is 
concerned with the more interesting non-atomic case and in particular 
shows that the property that the exponent is zero or positive is 
topologically invariant (we shall show below that the negative 
Lyapunov exponent case always corresponds to an atomic measure). 

The integrability of \( \log |Df| \) means that our definition of 
 Lyapunov exponents, commonly used in the one-dimensional context, 
 agrees with the more classical definition in terms of the limit of 
 the rate of growth of the derivative. Indeed, a standard application 
 of Birkhoff's ergodic theorem (which relies on the integrability 
 property) gives 
 \[ 
\lim_{n\to\infty} \frac{1}{n}\log |Df^{n}(x)| = \int\log |Df| d\mu 
=\lambda (\mu) 
\quad\text{ for } \mu \text{-a.e. } x. 
 \]
This pointwise definition can be generalised to 
the so-called upper and lower Lyapunov exponents 
\[  \lambda_{-}(x) := 
 \liminf_{n\to \infty} \frac1n \log |Df^n(x)|  
 \quad\text{and}\quad
 \lambda_{+}(x) := 
 \limsup_{n\to \infty} \frac1n \log |Df^n(x)|  
 \]
These quantities are defined at every point and a natural 
generalisation of the question answered above is whether the signs of 
these upper 
and lower Lyapunov exponents are topological invariants.
It was shown in \cite{NP} in the unimodal setting, that the positivity
of the lower Lyapunov exponent along the critical orbit (the
Collet-Eckmann condition) is preserved under topological conjugacy.
This result does not hold for multimodal maps, see \cite{PRS},
although it does generalise under additional recurrence conditions 
on the critical orbits \cite{LW}.
In \cite{PRS} it is also shown that in the context of rational
maps on the Riemann sphere, the property that the Lyapunov exponents of 
all invariant measures are uniformly positive is preserved under 
topological conjugacy. It is not known whether this extends to
$C^2$ interval maps.

If $f$ is unimodal and Collet-Eckmann, then every point has a positive 
upper Lyapunov exponent \cite{NowSan}. 
As the Collet-Eckmann condition is preserved
under conjugacy, the sign of upper pointwise Lyapunov exponent
is preserved under conjugacy for Collet-Eckmann maps. However we show 
that at least for lower Lyapunov exponents this is false in general. 

\begin{proposition}\label{lowlyap}
There exist unimodal maps with points for which
the sign of the lower pointwise Lyapunov exponent
is not preserved under topological conjugacy.
This is not restricted to orbits asymptotic to neutrally
attracting or neutrally repelling periodic orbits.
\end{proposition}

In \cite{PRS} this result was proved for bimodal maps;
their argument would not apply to the unimodal case, but shows that 
the lower pointwise Lyapunov exponent need not be preserved under 
a quasi-symmetric conjugacy. 

We make the following conjecture:

\begin{conjecture}\label{upplyap}
Topological conjugacy preserves the
sign of the upper pointwise Lyapunov exponents of all points
that are not attracted to a periodic orbit.
\end{conjecture}

It is immediate from the ergodic theorem that for every 
invariant measure $\mu$, there are points $x$ such that the
Lyapunov exponent $\lambda(\mu)$ coincides with the pointwise
Lyapunov exponent $\lambda(x)$.
(We write $\lambda(x)$ if $\lambda_+(x) = \lambda_-(x)$.)
However, there are instances
where a pointwise Lyapunov exponent 
is different from 
the Lyapunov exponent of all invariant measures.
This is shown in Proposition~\ref{lyap=0}.

\section{Proof of Theorem 
\ref{sign_lya_meas}}

\begin{lemma}\label{equidistribution}
If $\mu \in {\mathcal M}$ and \( \lambda(\mu)<0 \) then \( \mu \) is the 
Dirac-\( \delta \) measure equidistributed on an attracting periodic orbit. 
 \end{lemma}

\begin{proof} This follows from Przytycki's result \cite{Prz93} which states
that for $\mu \in {\mathcal M}$, $\lambda(x) \geq 0$ for $\mu$-a.e. $x$.
\end{proof}

For the proof of Theorem~\ref{sign_lya_meas}, we need a construction 
developed by Hofbauer \cite{Hofbauer}, called 
{\em canonical Markov extension}.
This Markov system is $(\hat I, \hat f)$, where $X$ is a disjoint union of 
closed intervals. 
Let $\P = \P_0 = \{ \xi_0 , \dots , \xi_r \}$ be the partition of
$I$ into the monotonicity intervals of $f$. Also write 
$\P_n = \bigvee_{i = 0}^{n-1} f^{-i}(\P_0)$, and $\P_n[x]$ is the element
of $\P_n$ containing $x$.
We will construct $X$ inductively.
\begin{itemize}
\item The {\em base} $B := I$ belongs to $\hat I$.
\item If $D \in \hat I$, let $E = \overline{ f(D \cap \xi_i) }$.
If the interval $E$ is equal to some already existing $D' \in \hat I$, then
define $\hat f(x,D) = (f(x), D')$. Otherwise, add $E$ disjointly to 
$\hat I$ and let $\hat f(x,D) = (f(x), E)$.
(Note that if $x \in \partial \xi_i$, then use 
$f(x) = \lim_{\xi_i \owns y \to x} T(y)$ to defined $f$ on $\partial \xi$.)
\end{itemize}
The system $(\hat I, \hat f)$ is Markov in the sense that the of any component
$D$ of $\hat I$ equals some union of components of $\hat I$. 
If we define the projection by
$\pi(x,D) = x$, then $\hat f \circ \pi = \pi \circ f$.
Due to the Markov property, the following is true.
\[
f^n(\P_n[x]) = D \in \hat I
\mbox{ if and only if } \hat f^n(\pi^{-1}(x) \cap B) \in D.
\]
If $\mu$ is $f$-invariant, then we can construct a measure $\hat \mu$ as follows:
Let $\hat \mu_0$ be the measure $\nu$ lifted to the base $B$
and set $\hat \mu_n = \frac{1}{n+1} \sum_{i=0}^{n} \hat \mu_0 \circ \hat
f^{-i}$. Clearly $\mu = \hat \mu_n \circ \pi^{-1}$ for each $n$.
As was shown in \cite{K1}, $\hat \mu_n$ converges vaguely
(i.e. on compact sets) to a limit measure, say $\hat \mu$.
If $\mu$ is ergodic, then $\hat \mu$ is either a probability measure 
on $\hat I$, in which case we call $\mu$ {\em liftable},
or $\hat \mu(D) = 0$ for all $D \in \hat I$.

\medskip

Let us say that an $n$-periodic point $p$ with 
{\em multiplier} $|Df^n(p)| \leq 1$ is 
{\em essential} if it is (one-sided) attracting and there 
exists $p' \in \orb(p)$ and a critical or boundary
point $c$ such that $f^i((c,p')) \cap {\mathcal C} = \emptyset$
for all $i \geq 0$.
This applies for example to periodic points of multimodal maps with 
negative Schwarzian derivative.

\begin{proposition}\label{liftable}
Let $f$ be a $C^3$ multimodal interval map with non-flat critical points. 
Let $\mu$ be an ergodic invariant probability measure such that
$\mu(p) = 0$ for each periodic point $p$ that is
(i) inessential with multiplier $\leq 1$, or
(ii) belongs to the boundary of the basin of another 
periodic point. 
Then $\mu$ is liftable if and only if $\mu$ has
a positive Lyapunov exponent.
\end{proposition}

This result was proved in \cite{BK}  for unimodal maps with negative Schwarzian
derivative. Here we give the 
details for the multimodal case, although the idea of proof is
the same.

\begin{proof}
The ``if'' part is proved in \cite{K1} using a construction 
from \cite{ledrappier}, except that \cite{K1} does not cover the case 
of atomic measures. 
So let us assume that $\mu$ is the equidistribution of a hyperbolic 
repelling periodic orbit $\orb(p)$, and that $p$ is not a boundary point 
of the basin of a periodic attractor.
Let $N = 2\times$ the period of $p$, so $f^N$ is orientation
preserving in a neighbourhood of $p$.
Let $Z_k$ be the largest neighbourhood of $p$ on which $f^{kN}$ is 
monotone. Write $Z_k = (a_k, b_k)$.\footnote{If $p$ is 
a boundary point of the interval $I$, then we have to adjust this argument to 
one-sided neighbourhoods $(a_k , p]$ or $[p,b_k)$.} 
Because $p$ is not a boundary point of the basin of a periodic attractor, 
$a_k$ and $b_k$ are precritical points. More precisely,
there are $n_a, n_b < N$
such that $f^{n_a}(a_1), f^{n_b}(b_1) \in {\mathcal C}$,
and because $f^N(Z_{k+1}) = Z_k$, $a_k \in f^{-n_a-(k-1)N}({\mathcal C})$
and $b_k \in f^{-n_b-(k-1)N}({\mathcal C})$.
It follows that if $p_0 = \pi^{-1}(p) \cap B$, then
$\hat f^{ik}(p_0) \in \pi^{-1}(p) \cap K_N$, where $K_N$ is the compact
part of $\hat I$ consisting of all components $D$ 
that can be reach by a path 
$B \to D_1 \to \dots \to D$ of length $\leq N$.
Clearly $\pi^{-1}(p) \cap K_N$ is finite, so it contains a $N$-periodic 
point $\hat p$. Because the lift $\hat \mu$ of $\mu$ is unique (see \cite{K1}),
$\hat\mu$ must be equal to the equidistribution on $\orb(\hat p)$.

For the ``only if'' part let us start proving that the equidistribution 
on a stable or neutral periodic orbit is non-liftable. 
Since such an orbit is essential, there is a point $p$ 
in this orbit and a critical (or boundary) point $c$ such that
$f^n( (c,p) ) \cap {\mathcal C} = \emptyset$ for all $n \geq 0$.
Assume that $p \in \xi_k \in \P$, and let
$\hat p = \pi^{-1}(p) \cap B$ be the lift of $p$ to the base $B$ of 
the Markov extension. Then $\hat f(\hat p)$
belongs to a successor $D := f(\xi_k)$ of $B$, and $f(c) \in \partial D$.
But since $f^n((c,p))$ never intersects a critical point,
each $\hat f^n(\hat p)$ belongs a different component of $\hat I$.
Therefore, the limit of the measures $\hat \mu_n$ is not liftable.

Let us assume that $\mu$ is
liftable, $\hat \mu$ being the lifted measure. We will show that
$\lambda(\mu) > 0$.
Let $D \in \hat I$ be such that $\hat \mu(D) > 0$ and let
$J$ be an interval, compactly contained in  $D$, such that $\hat \mu(J) > 0$.
Since $\mu$ is not the equidistribution on the orbit of a stable or neutral
periodic point $p$, $\pi(\overline J)$ can be chosen disjoint from $\orb(p)$.
Moreover we can chose $J$ such that $orb(\partial J) \cap J = \emptyset$.
Let $\hat F:J \to J$ be the first return map to $J$.
By our conditions on $J$ each branch $\hat F:J_i \to J$ of $\hat F$ is 
onto, and by the Markov property
of $\hat f$, $\hat F|_{J_i}$ is extendible monotonically to a branch 
that covers $D$.
Clearly each branch of $\hat F$, say $\hat F|_{J_i} = \hat f^s|_{J_i}$, 
contains an $s$-periodic point $q$. 
Due to a result by Martens, de Melo and van Strien \cite{MMS} and also
\cite[Theorem IV B']{dMvS},
there exists $\eps > 0$ such that the $|(\hat f^s)'(q)| > 1+\eps$,
independently of the branch. If $J$ is sufficiently small, the 
Koebe Principle \cite[Section IV.1]{dMvS}yields that
$|\hat F'(x)| > 1+\frac{\eps}{2}$ for all $x \in J$.
The Koebe Principle holds for maps with negative Schwarzian derivative,
but the work of Kozlovski \cite{Koz} and in the multimodal setting
van Strien \& Vargas \cite{SV} implies that the branches of $F$ have 
negative Schwarzian derivative if $f$ is $C^3$ and $J$ sufficiently small.
It is at this moment that we use the $C^3$ assumption.
It follows from the thesis of Mike Todd \cite{todd} that
a $C^{2 + \eps}$ assumption suffices for unimodal maps.

Clearly $\frac{\hat \mu}{\hat \mu(J)}$ is an
$\hat F$-invariant probability measure on $J$.
Let $J_i$, $i \in \N$, be the branch-domains of $\hat F$, and
let $s_i$ be such
that $\hat F|_{J_i} = \hat f^{s_i}|_{J_i}$.
Since we can write $\hat \mu$ 
as 
\[
\int \phi \ d\hat\mu = \sum_i
\sum_{j=0}^{s_i-1} \int_{J_i} \ph \circ \hat f^j d\hat \mu,
\]
we get
\begin{eqnarray*}
\lambda(\hat \mu) &=& \sum_i \sum_{j=0}^{s_i-1}
\int_{J_i} \log\, |\hat f'| \circ \hat f^j d\hat \mu \\
&=& \sum_i \int_{J_i} \log\, |(\hat f^{s_i})'| d\hat \mu
= \sum_i \int_{J_i} \log |\hat F'| d\hat \mu \\
&\geq&  \sum_i \hat \mu(J_i) \log(1+\frac{\eps}{2}) = \hat \mu(J) \log(1+\frac{\eps}{2}) > 0.
\end{eqnarray*}
Because $f'(\pi(x)) = \hat f'(x)$ for all $x \in \hat I$,
this concludes the proof.
\end{proof}

Now we are ready to prove Theorem~\ref{sign_lya_meas}.

\begin{proof}
First assume that $\lambda(\mu) > 0$, and let $\hat \mu$ be its lift to 
the Markov extension.
Assume that $g:J \to J$ has lift $\hat g:\hat J \to \hat J$.
Define $\hat h:\hat I \to \hat J$ as $\hat h|(D_n \subset \hat I) =
h|(D_n \subset I)$. Then $\hat h^* \hat \mu$ is a $\hat g$-invariant 
probability measure with $h^*\mu = \hat h^*\hat\mu \circ \pi^{-1}$.
Thus $\hat h^*\hat \mu$ is the lift of $h^*\mu$. It follows that 
$h^*\mu$ is liftable and hence has a positive Lyapunov exponent.
(Here we should recall that $\mu$ and $h^*\mu$ are non-atomic, so they are 
not associated with any periodic orbit, essential or not.)

By Lemma~\ref{equidistribution}, negative Lyapunov exponents can
only occur for atomic measures. Hence the remaining case $\lambda(\mu) = 0$
is also preserved under conjugacy.
\end{proof}

\section{Pointwise Lyapunov exponents}

\begin{proof}[Proof of Proposition \ref{lowlyap}]
We give a counter-example based on the unimodal maps 
$f(x) = 4x(1-x)$ and $g(x) = \sin(\pi x)$, both having negative Schwarzian 
derivative. These maps are conjugate 
on the unit interval. Due to the well-known smooth 
conjugacy with the tent map, we have that $\lambda(x) = \log 2$
whenever the limit exists and $f^n(x) \neq 1$ for all $n \geq 1$.
However, the limit need not always exist.
Indeed, let $(n_k)$ be a superexponentially 
increasing integer sequence, and $y \in [0,1]$ a point such that
\begin{itemize}
\item $y_i \in [\frac12,1]$ for $y_i = f^i(y)$ and $0 < i < n_1$.
Assuming $n_1$ is large, 
this means that $y_i \approx p = \frac34$, the fixed point of $f$, 
and hence $Df^{n_1}(y) \approx 2^{n_1}$.
\item $y_{n_1}$ is close to $c$ such that $y_{n_1} \approx 1$ and
$y_i \in [0,\frac12]$ for $n_1+1 < i \leq 2.1n_1$.
This means that $y_i \approx 0$, the other fixed point, and since 
$f'(0) = 4$, we obtain that $|0-y_{n_1+2}| = \O(4^{-1.1n_1})$.
Consequently, $|c-y_{n_1}| =  \O(\sqrt{4^{-1.1n_1}}) = \O(2^{-1.1n_1})$, 
and hence $|Df^{1+n_1}(y)| = \O(2^{n_1} \cdot 2^{-1.1n_1}) = \O(2^{-0.1n_1})$,
whereas $|Df^{2.1n_1}(y)| = 
\O(2^{n_1} \cdot 2^{-0.1 n_1} \cdot 4^{1.1 n_1}) = \O( 2^{2.1n_1} )$.
\item Let $y_i \in [\frac12,1]$ for $2.1n_1 < i < n_2$.
Hence, we find $Df^{n_2}(y) =  \O(2^{n_2})$.
\item $y_{n_2}$ is close to $c$ (and hence $y_{1+n_2}$ close to $1$) 
such that $y_i \in [0,\frac12]$ for $n_2+1 < i \leq 2.1n_2$.
It follows that
$|Df^{1+n_2}(y)| = \O(2^{n_1} \cdot 2^{-1.1n_2} ) = \O(2^{-0.1 n_2})$,
and $|Df^{2.1n_1}(y)| = \O(2^{2.1n_2})$.
\end{itemize}
Continue in this fashion, and we find that the lower
Lyapunov exponent is
$\underline \lambda(y) = \liminf \frac1n \log Df^n(y) = -0.1 \log 2$
whereas the upper Lyapunov exponent
$\overline \lambda(y) = \limsup \frac1n \log Df^n(y) = \log 2$

Now we do the same for $g = h \circ f \circ h^{-1}$ and the corresponding 
$\tilde y = h(y)$, we have to deal with different
multipliers:
$|Dg(0)| = \pi < |Df(0)|$ and $\alpha := |Dg(\tilde p)| \approx 2.12 > 
|Df(p)|$ for $\tilde p = h(p)$.
We now get that
$|Dg^{1+n_k}(\tilde y)| = \O( (\frac{\alpha}{ {\pi}^{0.55} })^{n_k} )$
is still exponentially large, 
so in this case, $\underline \lambda(\tilde y) > 0$.
\end{proof}

\subsection*{Example:}

We want to compare the results in this paper to an example from \cite{B}.
In this example, two conjugate smooth unimodal maps $f_1$ and $f_2$
(in fact, $f_1$ is quadratic and $f_2$ is a sine function),
for which
\[
\inf_{\eps > 0} \ \lim_{n\to\infty} 
\frac1n \{ 0 \leq i < n \ : \ f_k^i(c) \in (p -\eps,p +\eps) \}
= 1
\]
for $k = 1,2$ and $p = p_k$ is the orientation reversing fixed 
point of $f_k$.
Yet $f_1$ has an {\em acip} (i.e. an absolutely continuous (w.r.t Lebesgue) 
invariant probability measure), and $f_2$ has not.
Clearly the Dirac measure $\delta_p$ is the only weak limit point of 
$(\frac1n \sum_{i=0}^{n-1} \delta_{f_k^i(c)})$ for $k = 1,2$.
Any non-liftable measure
belongs to the convex hull of weak accumulation points of
$(\frac1n \sum_{i=0}^{n-1} \delta_{f^i(c)})$, see \cite{K1}.
Consequently, $f_k$ has only liftable invariant measures, all of
which have positive Lyapunov exponents.
The acip of $f_1$ does not transform under $h^*$ to an acip of $f_2$,
and in fact, there is not a single $f$-invariant
measure $\mu$ such that $h^*\mu$ is absolutely continuous.

A result by Keller \cite{K2} implies that for $k = 2$,
$\delta_p$ is the only weak limit 
point of  $(\frac1n \sum_{i=0}^{n-1} \delta_{f_2^i(x)})$
for Lebesgue-a.e. $x$. 
Recall that a {\em physical measure} $\mu$ is defined
by the fact that for every continuous observable $\phi:[0,1] \to \R$,
\begin{equation}\label{physical}
\mu(\phi) := \int \phi \ d\mu = 
\lim_{n\to\infty} \frac1n \sum_{i=0}^{n-1} \phi \circ f^i(x)
\mbox{ Lebesgue-a.e.}
\end{equation} 
Therefore $\delta_p$ is the physical measure of $f_2$.
However $\liminf_n \frac1n \log |Df^n(x)| = 0$ Lebesgue-a.e., because
otherwise there would be an acip by \cite{K2}.
This shows that it is important in \eqref{physical} to have continuous,
not just $L^1$, observables.

Since $f_k$ is not Collet-Eckmann,
$\inf\{ \lambda(\mu) : \mu \mbox{ is $f_k$-invariant} \} = 0$
for $k = 1,2$.
Therefore the infimum of Lyapunov exponents is not attained. This is in 
contrast to the Lyapunov exponent of invariant measures supported on hyperbolic sets, see \cite{CLR}.
The below results shows that the spectrum of pointwise Lyapunov exponents
can be strictly larger than the spectrum of Lyapunov exponents of 
measures.

\begin{proposition}\label{lyap=0}
There exists a unimodal map $f$ such that $\lambda(\mu) > 0$
for every $\mu \in {\mathcal M}$, but there is a point $x$
whose Lyapunov exponent exists (as a limit) and equals $0$.
\end{proposition}

\begin{proof}
We start by introducing some notation for unimodal maps.
A point $z < c$ is called a {\em precritical point closest to $c$} 
if $f^S(z) = c$ 
for some iterate $S$ and $f^i(c,z) \not\owns c$ for $0 \leq i \leq S$. 
There is an increasing sequence $(z_k)$ of such precritical points,
starting with $z_0 \in f^{-1}(c)$. The corresponding iterates $S_k$
such that $f^{S_k}(z_k) = 0$ are called {\em cutting times}. Clearly
$S_0 = 1$ and $S_k > S_{k-1}$ for each $k \geq 1$.
Let $U_k = (z_k, z_{k+1})$ and $\hat U_k  = (\hat z_{k+1}, \hat z_k)$ 
the interval at the other side of $c$ such that $f(U_k) = f(\hat U_k)$.
Note that the intervals $(z_{k-1},c)$ and $(c,\hat z_{k-1})$
are the largest intervals adjacent to $c$ on which $f^{S_k}$ is 
a diffeomorphism.

If $f$ has no periodic attractor, then $z_k \to c$.
If there is a $b$-periodic attractor and $B$ is the component of its
basin of attraction containing $c$, 
then $z_k \to \partial B$.
In fact, if $f$ has a neutrally attracting periodic orbit 
(at a saddle node bifurcation), then $\partial B$ 
contains a point of this orbit, and
$z_k$ converges to $\partial B$ in a polynomial way
(the precise rate of convergence depends on degeneracy of the 
neutral periodic orbit).

In \cite{B} this phenomenon is exploited by creating a cascade of almost 
saddle node bifurcations; there is an infinite sequence of integers $b_n$
and a map $f$ created as the limit of a sequence of maps $f_n$, where
$f_n$ has a $b_n$-periodic orbit at a saddle node bifurcation.
While perturbing $f_n$ to $f_{n+1}$, the geometric properties of the 
sequence $(z_k)$ is preserved to some extend.
In the example constructed in \cite{B}, the geometry of
$(z_k)$ is such that $|z_k-z_{k+1}|$ decreases polynomially
for values of $k$ associated to almost saddle node bifurcations, and
$|z_k-z_{k+1}|$ decreases exponentially for other values
of $k$. One can construct examples where the first behaviour dominates
such that the following properties hold:
\begin{enumerate}
\item $1 \leq S_k - S_{k-1} \leq 2$ for all $k \geq 1$; hence
$k < S_k \leq 2k$ for $k \geq 1$.
\item The distances $|f^{S_k}(c) - f^{S_k}(z_{k+1})|$,
$|f^{S_k}(z_{k+1}) - f^{S_k}(z_{k})|$ and
$|f^{S_k}(z_{k}) - f^{S_k}(z_{k-1})|$ are bounded away from $0$, uniformly
in $k$. Using the Koebe Principle \cite{dMvS}, we conclude that the distortion
of $f^{S_k}| U_k$ and  $f^{S_k}| \hat U_k$ is uniformly bounded.
\item 
$\lim_k \frac1k \log |z_k - z_{k+1}|^{-1} = 0$.
\end{enumerate}
Construct the induced map $F$ by
$F|U_k \cup \hat U_k = f^{S_k}$. It is easy to verify from
property (1) that $F(U_k) = F(\hat U_k) = (z_0,c)$, $(z_1,c)$, $(c,\hat z_0)$ 
or $(c,\hat z_1)$.
Hence $F$ is a Markov map.
For any $x$, write $\chi_n(x) = k$ if $F^n(x) \in U_k \cap \hat U_k$.
Also, let $t_n = \sum_{i=0}^{n-1} S_{\chi_i(x)}$, so
$x_n := F^n(x) = f^{t_n}(x)$.
Because of the Markov properties of $F$, there are points $x$ such that
$\chi_n(x) \to \infty$ so slowly that $\frac{t_{n+1}-t_n}{t_n} \to 0$.
Therefore
\begin{eqnarray*}
\frac{1}{t_n} \log |Df^{t_n}(x)| 
&=& \frac{1}{t_n} \log \prod_{i=0}^{n-1} |Df^{S_{\chi_i}}(x_i)| \\
&=& \frac{ \sum_{i=0}^{n-1} \log  |Df^{S_{\chi_i}}(x_i)| } 
{\sum_{i=0}^{n-1} S_{\chi_i}  }\\
&\sim& \frac{ K \ \sum_{i=0}^{n-1}\log |z_{\chi_i} - z_{\chi_i+1}|^{-1} }
{\sum_{i=0}^{n-1} \chi_i  } \to 0,
\end{eqnarray*}
by property (3).
Here $K$ depends only on the image-length and distortion of the branches 
of $F$, which are uniform by property (2).
Finally, for intermediate values of $t$, i.e. $t_n \leq t < t_{n+1}$,
we have 
\[
L^{ t - t_{n+1} } \ |Df^{t_{n+1}}(x)| \leq |Df^t(x)| 
\leq L^{t-t_{n} } \ |Df^{t_n}(x)|
\] 
for $L = \sup |Df| < \infty$. By the assumption that  
$\frac{t_{n+1}-t_n}{t_n} \to 0$, we obtain
$\lim_t \frac{1}{t} \log |Df^{t}(x)| = 0$ as well.
This concludes the proof.
\end{proof}

\end{document}